\documentclass[12pt]{article}
\usepackage{amssymb,amsmath,epsf}

\newcommand{\figone}{

\begin{flushright}
\unitlength=1mm
\begin{picture}(140,71)(5,-71)
\thicklines

\put(7,-3)
{\makebox(0,0)[l]{
\it {\bf Figure 1:} A smooth projective curve $Z$ of genus $g$.}}

\put(7,-9){
\parbox[t]{120cm}{
\small \footnotesize
A rational curve is an analogue of a closed 
interval with two boundary points. \\
An elliptic curve is an analogue of a circle. 
Higher genus curves correspond to graphs.}}

\put(8,-25){\makebox(0,0){$g=0$}}
\put(55,-25){\makebox(0,0){$g=1$}}
\put(110,-25){\makebox(0,0){$g>1$}}

\put(-2,-35){\line(1,0){20}}
\put(-2,-35){\circle*{1.3}}
\put(18,-35){\circle*{1.3}}

\put(55,-35){\circle{12}}

\put(87,-29){\line(1,0){29}}
\put(132,-29){\line(-1,0){5}}
\put(87,-41){\line(1,0){29}}
\put(132,-41){\line(-1,0){5}}

\multiput(119,-29)(2.5,0){3}{\circle*{0.5}}
\multiput(119,-41)(2.5,0){3}{\circle*{0.5}}

\put(87,-35){\oval(12,12)[l]}
\put(132,-35){\oval(12,12)[r]}

\multiput(89,-41)(12,0){3}{\line(0,1){12}}
\multiput(89,-29)(12,0){3}{\circle*{1.2}}
\multiput(89,-41)(12,0){3}{\circle*{1.2}}

\put(130,-41){\line(0,1){12}}
\multiput(130,-41)(0,12){2}{\circle*{1.2}}

\put(-14,-48){
\parbox[t]{43mm}{
\footnotesize
$Z$ is polar orientable (with polar boundary):
there exists a 1-form with two simple poles and without zeros,
}}

\put(8,-65){\makebox(0,0){$ hp_0(Z)=1\,,$}}
\put(8,-70){\makebox(0,0){$ hp_1(Z)=0\,,$}}

\put(35,-48){
\parbox[t]{40mm}{
\footnotesize
$Z$ is polar orientable (without a polar boundary):
there exists a 1-form without zeros or poles,
}}

\put(55,-65){\makebox(0,0){$ hp_0(Z)=1\,,$}}
\put(55,-70){\makebox(0,0){$ hp_1(Z)=1\,,$}}

\put(88,-49){
\parbox[t]{42mm}{
\footnotesize
$Z$ is not polar orientable:
a generic holomorphic 1-form has $2g-2$ zeros,
}}

\put(110,-65){\makebox(0,0){$ hp_0(Z)=1\,,$}}
\put(110,-70){\makebox(0,0){$ hp_1(Z)=g\,.$}}

\end{picture}
\end{flushright}
}

\newcommand{\figtwo}{

\begin{flushright}
\unitlength=1mm
\begin{picture}(140,47)(5,-47)
\thicklines

\put(7,-3)
{\makebox(0,0)[l]{
\it {\bf Figure 2:} A smooth projective curve without a point, 
$Z\smallsetminus \{P\}$.}}

\put(8,-15){\makebox(0,0){$g=0$}}
\put(55,-15){\makebox(0,0){$g=1$}}
\put(110,-15){\makebox(0,0){$g>1$}}

\put(-2,-25){\line(1,0){9}}
\put(18,-25){\line(-1,0){9}}
\put(-2,-25){\circle*{1.3}}
\put(18,-25){\circle*{1.3}}

\put(48.42,-31.66){\epsffile{S1-P.eps}}

\put(87,-19){\line(1,0){7}}
\put(96,-19){\line(1,0){20}}
\put(132,-19){\line(-1,0){5}}
\put(87,-31){\line(1,0){29}}
\put(132,-31){\line(-1,0){5}}

\multiput(119,-19)(2.5,0){3}{\circle*{0.5}}
\multiput(119,-31)(2.5,0){3}{\circle*{0.5}}

\put(87,-25){\oval(12,12)[l]}
\put(132,-25){\oval(12,12)[r]}

\multiput(89,-31)(12,0){3}{\line(0,1){12}}
\multiput(89,-19)(12,0){3}{\circle*{1.2}}
\multiput(89,-31)(12,0){3}{\circle*{1.2}}

\put(130,-31){\line(0,1){12}}
\multiput(130,-31)(0,12){2}{\circle*{1.2}}

\put(8,-40){\makebox(0,0){$ hp_0(Z\smallsetminus \{P\})=2\,,$}}
\put(8,-45){\makebox(0,0){$ hp_1(Z\smallsetminus \{P\})=0\,,$}}

\put(55,-40){\makebox(0,0){$ hp_0(Z\smallsetminus \{P\})=1\,,$}}
\put(55,-45){\makebox(0,0){$ hp_1(Z\smallsetminus \{P\})=0\,,$}}

\put(90,-40){\makebox(0,0)[l]{$ hp_0(Z\smallsetminus \{P\})=1\,,$}}
\put(90,-45){\makebox(0,0)[l]{$ hp_1(Z\smallsetminus \{P\})=g-1\,.$}}

\end{picture}
\end{flushright}
}

\newcommand{\figthree}{

\begin{center}
\unitlength=1mm
\begin{picture}(140,83)(5,-83)
\thicklines

\put(70,-3)
{\makebox(0,0)[t]{
\it {\bf Figure 3:} A smooth projective curve without two points, 
$Z\smallsetminus \{P,Q\}$.}}

\put(0,-20){\makebox(0,0)[l]{$g=0$}}
\put(0,-35){\makebox(0,0)[l]{$g=1$}}
\put(0,-55){\makebox(0,0)[l]{$g>1\,,~P,Q~ \mbox{\small are generic}$}}
\put(0,-75){\makebox(0,0)[l]{$g>1\,,~P+Q~ \mbox{\small is special}$}}

\put(50,-20){\line(1,0){6}}
\put(70,-20){\line(-1,0){6}}
\put(58,-20){\line(1,0){4}}
\put(50,-20){\circle*{1.3}}
\put(70,-20){\circle*{1.3}}

\put(53.5,-41.5){\epsffile{S1-PQ.eps}}

\put(50,-49){\line(1,0){7}}
\put(59,-49){\line(1,0){10}}
\put(71,-49){\line(1,0){8}}
\put(95,-49){\line(-1,0){5}}
\put(50,-61){\line(1,0){29}}
\put(95,-61){\line(-1,0){5}}

\multiput(82,-49)(2.5,0){3}{\circle*{0.5}}
\multiput(82,-61)(2.5,0){3}{\circle*{0.5}}

\put(50,-55){\oval(12,12)[l]}
\put(95,-55){\oval(12,12)[r]}

\multiput(52,-61)(12,0){3}{\line(0,1){12}}
\multiput(52,-49)(12,0){3}{\circle*{1.2}}
\multiput(52,-61)(12,0){3}{\circle*{1.2}}

\put(93,-61){\line(0,1){12}}
\multiput(93,-61)(0,12){2}{\circle*{1.2}}

\put(50,-69){\line(1,0){7}}
\put(59,-69){\line(1,0){20}}
\put(95,-69){\line(-1,0){5}}
\put(50,-81){\line(1,0){7}}
\put(59,-81){\line(1,0){20}}
\put(95,-81){\line(-1,0){5}}

\multiput(82,-69)(2.5,0){3}{\circle*{0.5}}
\multiput(82,-81)(2.5,0){3}{\circle*{0.5}}

\put(50,-75){\oval(12,12)[l]}
\put(95,-75){\oval(12,12)[r]}

\multiput(52,-81)(12,0){3}{\line(0,1){12}}
\multiput(52,-69)(12,0){3}{\circle*{1.2}}
\multiput(52,-81)(12,0){3}{\circle*{1.2}}

\put(93,-81){\line(0,1){12}}
\multiput(93,-81)(0,12){2}{\circle*{1.2}}

\put(105,-19){\makebox(0,0)[bl]{$ hp_0(Z\smallsetminus \{P,Q\})=3\,,$}}
\put(105,-24){\makebox(0,0)[bl]{$ hp_1(Z\smallsetminus \{P,Q\})=0\,,$}}

\put(105,-34){\makebox(0,0)[bl]{$ hp_0(Z\smallsetminus \{P,Q\})=2\,,$}}
\put(105,-39){\makebox(0,0)[bl]{$ hp_1(Z\smallsetminus \{P,Q\})=0\,,$}}

\put(105,-54){\makebox(0,0)[bl]{$ hp_0(Z\smallsetminus \{P,Q\})=1\,,$}}
\put(105,-59){\makebox(0,0)[bl]{$ hp_1(Z\smallsetminus \{P,Q\})=g-2\,,$}}

\put(105,-74){\makebox(0,0)[bl]{$ hp_0(Z\smallsetminus \{P,Q\})=2\,,$}}
\put(105,-79){\makebox(0,0)[bl]{$ hp_1(Z\smallsetminus \{P,Q\})=g-1\,.$}}

\end{picture}
\end{center}
}

\newcommand{\mysection}[1]{\section{\large\bf #1}}

\newtheorem{defn}[subsection]{Definition}
\newtheorem{prop}[subsection]{Proposition}
\newtheorem{theo}[subsection]{Theorem}

\newenvironment{ssect}[1]{\smallskip\noindent%
\refstepcounter{subsection}%
{\bf \thesubsection~#1}\hspace{-1mm}}
{\smallskip}

\newenvironment{rem}{\smallskip\noindent%
\refstepcounter{subsection}%
{\bf \thesubsection}~~{\sc Remark.}\hspace{-1mm}}
{\smallskip}

\newenvironment{ex}{\smallskip\noindent%
\refstepcounter{subsection}%
{\bf \thesubsection}~~{\sc Example.}}
{\smallskip}

\newenvironment{ackn}{\medskip \noindent \small
{\sl Acknowledgments.}}{\bigskip}

\newenvironment{smallbibl}[1]
{\small
\begin{thebibliography}{#1}}
{\end{thebibliography}}

\newcommand{\res}{{\rm res}}

\newcommand{\dvsr}{{\rm div}}
\newcommand{\lk}{\ell k_{polar}}
\newcommand{\C}{{\mathbb C}}
\newcommand{\R}{{\mathbb R}}

\renewcommand{\P}{{\mathbb P}}

\newcommand{\Cdot}{\,\mbox{{\large$\cdot$}}\,}
\newcommand{\lbar}[1]{\,\overline{\!#1}}
\newcommand{\tpi}{2\pi i}

\begin{document}

\title{
\vspace{-2.4cm}
\begin{flushright}
{\normalsize{\it Philos. Trans. Roy. Soc. London Ser. A,} 2001}\\
\vspace{-0.2cm}
{\normalsize ITEP-TH-05/01}\\
\end{flushright}
\vspace{0.7cm}
{\Large\sc Polar Homology and Holomorphic Bundles}
}
\author{Boris Khesin\thanks{Department of Mathematics,
University of Toronto,
Toronto, ON M5S 3G3, Canada;
e-mail: {\tt khesin@math.toronto.edu}}~
and Alexei Rosly\thanks{Institute of Theoretical and
Experimental Physics, B.Cheremushkinskaya 25, 117259 Moscow,
\mbox{Russia}; e-mail: {\tt rosly@heron.itep.ru}}}

\date{December 20, 2000}

\maketitle
\begin{abstract}
We describe polar homology groups for complex  manifolds.
The polar $k$-chains are subvarieties of complex
dimension $k$ with meromorphic forms on them,
while the boundary operator is defined by
taking the polar divisor and the Poincar\'e
residue on it. The polar homology groups may be regarded as
holomorphic analogues of the homology groups in topology.
We also describe the polar homology groups
for quasi-projective one-dimensional varieties 
(affine curves). These groups obey
the Mayer--Vietoris property.
A complex counterpart
of the Gauss linking number of two curves in a three-fold
and various gauge-theoretic aspects of the above correspondence
are also discussed.

\end{abstract}


\mysection{Introduction} \label{Int}

In this paper we  describe certain homology groups,
for complex projective and one-dimensional quasi-projective manifolds. 
These {\it polar homology~} groups can be regarded  as a complex 
geometric counterpart  of singular homology groups in topology. 

The essence of the ``polar homology" theory described below
is presented in the following ``complexification dictionary":
$$
\begin{array}{rcl}
\text{a~real~manifold} &\leftrightarrow & 
\text{a~complex~manifold}\\
\text{an~orientation~of~the~manifold} &\leftrightarrow &
\text{a~meromorphic~volume~form~on~the~manifold}\\
\text{manifold's~boundary} &\leftrightarrow & 
\text{form's~divisor~of~poles}\\
\text{induced~orientation~of~the~boundary} &\leftrightarrow & 
\text{residue~of~the~meromorphic~form}\\
\text{ open~manifold's~infinity} &\leftrightarrow & 
\text{form's~divisor~of~ zeros}\\
\text{Stokes~formula} &\leftrightarrow & 
\text{Cauchy~formula}\\
\text{singular~homology} &\leftrightarrow & 
\text{polar~homology}\\
\end{array}
$$

\medskip

In short, polar $k$-chains in a complex projective manifold 
are  linear combinations of $k$-dimensional complex submanifolds 
with meromorphic
closed $k$-forms on them. The boundary operator sends such a pair \it (complex 
submanifold,  meromorphic form)\rm~ to the pair \it(form's divisor of poles,
form's residue at  the divisor),\rm~ that is, to a $(k-1)$-chain in
the same ambient manifold. The square of the boundary operator is
zero, and the polar homology groups  are defined as the quotients of
polar cycles over polar boundaries, see Section \ref{Formal}.

While the form's divisor of poles on a complex manifold is 
an analogue of  the boundary of a real manifold, the form's 
divisor of zeros can be related to the ``infinity'' of  a real
manifold, if the latter is non-compact
(see  Section \ref{QP}). 

\smallskip

This parallelism between topology and algebraic geometry extends to 
various gauge-theoretic notions and facts. In particular,
we discuss below several problems related to
the correspondence of  flat and holomorphic bundles.
Some features of this correspondence are  
also present in the papers \cite{Arn, FK, arhar, DT, Kh, T, KR}.
Note that the gauge theory  related 
to a version of the Chern--Simons functional on Calabi--Yau manifolds 
(see \cite{W}) was a motivation for the construction of these homology 
groups and of the relevant notion of the polar linking number (see
Section \ref{Gauge} and cf. \cite{FT, KR2}). 

\bigskip


\mysection{Polar homology of projective manifolds} \label{Proj}

We start with a heuristic motivation for polar homology and recall 
(following \cite{KR2}) the formal definition 
of the corresponding groups in the next section.

\begin{ssect}{A holomorphic analogue of orientation.} \label{InPair}
In order to see why a meromorphic or holomorphic form on a complex manifold
can be regarded as an analogue of orientation of a real manifold, we  
extend  the analogy between de Rham and Dolbeault cochains
(~~$d\leftrightarrow\bar\partial$~~) to an analogy at the level of the
corresponding {\it chain~} complexes.

Let $X$ be a compact complex manifold and $u$  be a smooth
$(0,k)$-form on it, $0\leqslant k\leqslant n=\dim\,X$.
We would like to treat such $(0,k)$-forms in the same manner as
ordinary $k$-forms on a smooth manifold, but in the framework
of complex geometry.  In particular, we have to be able to
integrate them  over $k$-dimensional {\it complex} submanifolds in $X$.
Recall that in the theory of differential forms, a form can be
integrated over a real submanifold provided
that the submanifold is endowed with an orientation. 
Thus, we need to find a holomorphic analogue of the orientation.

For a $k$-dimensional submanifold $W\subset X$ is equipped
with a holomorphic
$k$-form  $\omega$ one can consider the following integral
$$
\int_W \omega\wedge u \;  
$$
of the product of the $(k,0)$- and $(0,k)$-forms.
Therefore, here we are going to regard a top degree holomorphic 
form $\omega$ on a complex
manifold as an analogue of orientation. 

\end{ssect}

\begin{ssect}{The Cauchy--Stokes formula.}  \label{CSt}
More generally, if the form $\omega$
is allowed to have first order poles on a smooth hypersurface in $W$, the above
integral is still well-defined.
The new feature brought by the presence of poles of
$\omega$ manifests in the following relation.

Consider the integral $\int_W \omega\wedge u \;$ with 
a meromorphic $k$-form $\omega$
having first order poles on a smooth hypersurface $V\subset W$.
Let the smooth $(0,k)$-form $u$ on $X$ be
$\bar\partial\,$-exact, that is $u=\bar\partial\,v$ for some
$(0,k-1)$-form $v$ on $X$. Then
$$
\int_W \omega\wedge\bar\partial\,v
= 2\pi i\int_V\res\,\omega\wedge v \;.\label{CSf}
$$
We shall exploit this straightforward generalization of the Cauchy
formula as a comp\-lexi\-fied analogue of the Stokes theorem.

In the above formula  $\res\,\omega$ denotes a  $(k-1)$-form on $V$
which is the {\it Poincar\'e residue} of $\omega$.
Namely, the form $\omega$ can be locally expressed as
$\omega=\rho\wedge dz/z+\varepsilon,$
where $z=0$ is a local equation of $V$ in $W$ and $\rho$ (respectively, 
$\varepsilon$) is a holomorphic $(n-1)$-form (resp., $n$-form).
Then the restriction $\rho|_V$ is an unambiguously defined holomorphic 
$(n-1)$-form on $V$, and it is called the Poincar\'e residue
$\res\,\omega$ of the form $\omega$.  
\end{ssect}

\begin{ssect}{Boundary operator.}
The Cauchy--Stokes formula prompts  us to consider
the pair $(W,\omega)$ consisting of a $k$-dimensional submanifold
$W$ equipped with a meromorphic form $\omega$
(with first order poles on $V$) as an
analogue of a compact oriented submanifold with boundary. In the
polar homology theory  the
pairs $(W,\omega)$ will play the role of chains, while the boundary
operator will take the form $\partial\,(W,\omega)=\tpi(V,\res\,\omega)$.
Note, that in the situation under consideration, when the polar set
$V$ of the form $\omega$ is a smooth $(k-1)$-dimensional submanifold in a
smooth $k$-dimensional $W$, the induced ``orientation" on $V$ is given
by a holomorphic $(k-1)$-form $\res\,\omega$. This means that
$\partial\,(V,\res\,\omega)=0$, or the boundary of a boundary is zero.
The latter is the source of the identity $\partial^2=0$, which allows one to
define {\em polar homology} groups $HP_k$. 

\end{ssect}

\begin{ssect}{Pairing to smooth forms.} \label{InHom}
It is clear that the  polar homology groups of a complex
manifold $X$ should have a pairing to Dolbeault cohomology groups
$H^{0,k}_{\bar\partial}(X)$. Indeed, for a polar $k$-chain $(W,\omega)$
and any $(0,k)$-form $u$ such a pairing is given by the integral
$$
\langle(W,\omega)\,,u\rangle = \int_W \omega\wedge u \;.\label{pair}
$$
In other words, the polar chain $(W,\omega)$ defines a current on $X$ of
degree   $(n,n-k)$,
where $n=\dim X$. One can see that
this pairing descends to
(co)homology classes by virtue of the Cauchy--Stokes
formula:
$$
\langle(W,\,\omega)\,, \,{\bar\partial} v\rangle
=\langle \;\partial(W,\,\omega)\,, v\rangle.
$$
\end{ssect}

\begin{ex} \label{InCurv}
Now we are already able to find out the polar homology groups $HP_k$
of a complex projective curve $Z$.  In this (and in any)
case, all the 0-chains are cycles. Let $(P,a)$ and $(Q,b)$ be
two 0-cycles, where
$P,Q$ are points on $Z$ and $a,b\in\C\,$. They are polar
homologically equivalent iff $a=b$.  Indeed, $a=b$ is necessary and
sufficient for the existence of a meromorphic 1-form $\alpha$ on $Z$,
such that $\dvsr_\infty\alpha = P+Q$ and
$\res_P\,\alpha=\tpi\,a,\, \res_Q\,\alpha=-\tpi\,b$.
(The sum of all residues of a meromorphic differential on a
projective curve is zero by the Cauchy theorem.)
Then, we can write in
terms of polar chain complex (to be formally defined 
in the next section) that
$(P,a)-(Q,a)=\partial\,(Z,\alpha)$. Thus, $HP_0(Z)=\C\,$.

Polar 1-cycles correspond to all possible holomorphic
1-forms on $Z$. On the other hand, there are no 1-boundaries,
since there are no polar 2-chains in $Z$. Hence $HP_1(Z)\cong\C\,^g$, where
$g$ is the genus of the curve $Z$.
\end{ex}

\begin{ssect}{Polar intersections.} \label{X}
One can define a polar analogue  of the intersection number
in topology. For instance, let
$(X, \mu)$ be a complex manifold equipped with
a meromorphic volume form $\mu$ without zeros (its ``polar orientation").
Consider two polar cycles
$(A,\alpha)$ and $(B,\beta)$ of complimentary dimensions that
intersect transversely in $X$ (here $\alpha$ and $\beta$ are
volume forms, or ``polar orientations," on the corresponding submanifolds).
Then the polar  intersection number is defined by the formula
$$
        \langle(A,\alpha) \Cdot (B,\beta)\rangle =
    \sum_{P\in A\cap B}\frac{\alpha(P)\wedge\beta(P)}{\mu(P)} \;\,.
$$
At every intersection point $P$, the ratio in the right-hand-side
is the ``comparison" of the orientations of the polar cycles at that
point (the form $\alpha\wedge\beta$ at $P$)
with  the orientation of the ambient manifold (the form $\mu$ at $P$).
This is a straightforward analogue of the use of mutual orientation
of cycles in the definition of the topological intersection number.
Note, that in the polar  case the intersection number does not have
to be an integer.
(Rather, it is a holomorphic function of the ``parameters"
$(A,\alpha), (B,\beta)$ and $(X,\mu)$.)

Similarly, there is a polar analogue of the
intersection product of cycles when they intersect over
a manifold of positive dimension, given essentially by the same 
formula (see [KR2]). Furthermore, one can define a polar analogue 
of the linking number using the same philosophy of polar chains. 
We discuss polar linkings, which are very close in spirit to the 
polar intersections, in relation to the Chern--Simons theory at the 
end of the paper.
\end{ssect}

\begin{rem}
Most of the above  discussion extends to polar chains
$(A,\alpha)$  where the meromorphic $p$-form $\alpha$ is not
necessarily of top degree, that is $0\leqslant p\leqslant k$,
where $k=\dim_{\C} A$.
To define the boundary operator we have to restrict ourselves
to the meromorphic forms with logarithmic singularities.
The corresponding polar homology groups are enumerated 
by two indices $k$ and $p$ ($0\leqslant p\leqslant k$).
One can see that the Cauchy--Stokes formula extends to this case
as well, if we pair meromorphic $p$-forms $\omega$ on $W$ with smooth
$(k-p, p)$-forms on $X$.          
\end{rem}

\mysection{Definition of polar homology groups.}\label{Formal}

\begin{ssect}{Polar chains.}\label{smooth}
In this section we deal with complex projective varieties,
i.e., subvarieties of a complex projective space.
By a smooth projective variety we always understand a smooth and connected
one. For a smooth variety $M$, we denote by
$\Omega^p_M$ the sheaf of holomorphic $p$-forms on $M$.
The sheaf $\Omega^{\dim M}_M\,$ of forms of
the top degree on $M$ will sometimes be denoted by $K_M$.

         The space of polar $k$-chains for a complex
projective variety $X, \dim X=n,$
will be defined as a $\C\,$-vector space with certain generators and
relations.
\end{ssect}

\begin{defn} \label{PHdef}
The space of {\sl polar $k$-chains~} ${\cal C}_k(X)$ is
a vector space over $\C$ defined as the quotient
${\cal C}_k(X)=\hat{\cal C}_k(X)/{\cal R}_k$, where the vector space
$\hat{\cal C}_k(X)$ is freely generated by the
triples $(A,f,\alpha)$ described in (i),(ii),(iii) and
${\cal R}_k$ is defined as relations
{\rm (R1),(R2),(R3)} imposed on the triples.

(i) $A$ is a smooth complex projective variety, $\dim A=k$;

(ii) $f\!: A\to X$ is a holomorphic map of projective varieties;

(iii) $\alpha$ is a rational $k$-form on $A$
with first order poles on $V\subset A$, where $V$ is a normal crossing
divisor in $A$, i.e., $\alpha\in\Gamma(A,\Omega^k_A(V))$.

\noindent
The relations are:

{\rm (R1)} $\lambda (A, f, \alpha)=(A, f, \lambda\alpha)$

{\rm (R2)} $\sum_i(A_i,f_i,\alpha_i)=0$ provided that
 $\sum_if_{i*}\alpha_i\equiv 0$, where
$\dim f_i(A_i)=k$ for all $i$ and
the push-forwards $f_{i*}\alpha_i$ are considered on the smooth part
of $\cup_i f_i(A_i)$;\footnote{See, e.g., \cite{Gr} for the definition 
of the push-forward (or, {\it trace}) map on forms.}

{\rm (R3)} $(A,f,\alpha)=0$ if $\dim f(A)<k$.

\end{defn}

Note that by definition, ${\cal C}_k(X)=0$ for $k<0$ and 
$k>\dim X$.

\begin{rem}
The relation (R2) allows us, in particular, 
to deal with pairs instead of triples 
replacing a triple $(A,f,\alpha)$ by a pair
$(\hat A,\hat\alpha)$, where $\hat A=f(A)\subset X$, $\hat\alpha$
is defined only on
the smooth part of $\hat A$ and $\hat \alpha=f_*\alpha$ there.
Due to the relation (R2), such a pair $(\hat A,\hat\alpha)$ carries
precisely
the same information as $(A,f,\alpha)$.
(The only point to worry about
is that such pairs cannot be arbitrary. In fact, by the Hironaka
theorem on resolution of singularities, any subvariety $\hat A\subset X$
can be the image of some regular $A$, but the form $\hat\alpha$ on the
smooth part of $\hat A$ cannot be arbitrary.)

The same relation (R2) also represents additivity with respect to $\alpha$,
that is
$$
(A,f,\alpha_1)+(A,f,\alpha_2)=(A,f,\alpha_1+\alpha_2).
$$
Formally speaking, the right hand side makes sense only if
$\alpha_1+\alpha_2$ is an admissible form on $A$, that is if its polar
divisor $\dvsr_\infty(\alpha_1+\alpha_2)$ has normal crossings.
However, one can always replace $A$ with a variety $\tilde A$,
obtained from $A$ by a blow-up,
$\pi\!:\tilde A\to A$, in such a way
that $\pi^*(\alpha_1+\alpha_2)$ is admissible on $\tilde A$, i.e.,
$\dvsr_\infty(\alpha_1+\alpha_2)$ is already a normal crossing
divisor. (This is again the Hironaka theorem.) The (R2) says that
$(A,f,\alpha_1)+(A,f,\alpha_2)=(\tilde A,f\circ\pi,\pi^*(\alpha_1+\alpha_2))$.

\end{rem}

\begin{defn} \label{d-def}
The {\sl boundary operator~} $\partial: {\cal C}_k (X)\to{\cal C}_{k-1}(X)$
is defined by
$$
\partial(A,f,\alpha)=\tpi\sum_i(V_i, f_i, \res_{V_i}\,\alpha)
$$
(and by linearity),
where $V_i$ are the
components of the polar divisor of $\alpha$,
$\dvsr_\infty\alpha=\cup_iV_i$, and the maps $f_i=f|_{V_i}$
are restrictions of the map $f$ to each component of the divisor.

\end{defn}

\begin{theo}{\cite{KR2}}
The boundary operator $\partial$ is well defined, i.e.\ it is
compatible with the relations {\rm (R1),(R2),(R3)}. Moreover,
$~~\partial^2=0\;.$
\end{theo}

For the proof  we refer to \cite{KR2}. Note that having proved
compatibility, the relation $\partial^2=0$ becomes nearly evident. Indeed,
it suffices to prove it  for normal crossing divisors
of poles. In the latter case, the repeated residue 
at pairwise intersections differs by a sign according to
the order in which the residues are taken.
Thus the contributions to the repeated residue from different
components cancel out. 

\begin{defn}
For a smooth complex projective variety $X, \dim X=n$, the chain
complex
$$
0 \to {\cal C}_n(X) \stackrel{\!\!\partial}{\longrightarrow}
{\cal C}_{n-1}(X) \stackrel{\!\!\partial}{\longrightarrow} \dots
\stackrel{\!\!\partial}{\longrightarrow} {\cal C}_0(X)\to 0
$$
is called the {\sl polar chain complex} of $X$. Its homology groups,
$HP_k(X), k=0,\ldots,n$, are called the {\sl polar homology groups} of $X$.

\end{defn}

\begin{rem} \label{non-top}
As we mentioned before,  one can similarly define the  polar
homology groups $HP_{k,p}(M)$ for  the case of
$p$-forms on $k$-manifolds, i.e., for the forms of
not necessarily top degree, $p\leq k$.
Instead of meromorphic $k$-forms with poles of the first order we
have to restrict ourselves by $p$-forms with logarithmic
singularities, keeping the definition of the boundary operator $\partial$
intact. 
\end{rem}
\bigskip


\mysection{Polar Homology for Affine Curves}
\label{QP}

In the preceding section we introduced polar homology of projective
varieties. From the point of view of topological analogy (cf.\
Sect.~\ref{Int}), the projective varieties play the role of compact
spaces. It would be useful, of
course, to have also a consistent analogue of homology of arbitrary,
i.e.\ not necessarily compact, manifolds. It is natural to expect
that this latter role is played by Zariski open subsets in projective
varieties, that is by quasi-projective varieties. This is indeed the
case and the definition of polar homology can be extended to the
quasi-projective case, so that the polar homology groups 
obey certain natural properties expected from the topological
analogy. In particular, they obey the Mayer--Vietoris principle.

To simplify the exposition we shall describe here the
case of dimension one only, i.e., that of affine curves.

\smallskip

\begin{ssect}{}
Let $X$ be an affine curve and $\lbar{X}\supset X$ be its projective
closure. We shall define the
polar chains for the quasi-projective variety $X$ as a certain subset of polar
chains for $\lbar{X}$, but the result will depend only on $X$ and not
on the choice of $\lbar{X}$. Let us denote by $D$ the compactification
divisor, $D=\lbar{X}\smallsetminus X$. By differentials of the third
kind  on a complex curve we shall understand, as usual, meromorphic
1-forms which may have only first order poles.

\end{ssect}

\begin{defn} \label{QPdef}
The space ${\cal C}_0(X)$ is the vector space formed by 
complex linear combinations of points in $X$. It is a subspace in ${\cal
C}_0(\lbar{X})$. 

The vector space ${\cal C}_1(X)$ is defined as the subspace in
${\cal C}_1(\lbar{X})$ generated by the triples $(A,f,\alpha)$
where $A$ is a smooth projective curve, $f$ is a map
$f\!: A\to\lbar{X}$, and $\alpha$ is a differential of the third kind
on $A$ that vanishes at $f^{-1}(D)\subset A$.

\end{defn}

\begin{prop} \label{comp-indep}
The spaces ${\cal C}_k(X), k=0,1$,  
form a subcomplex in the polar chain complex
$({\cal C}_\bullet(\lbar{X}),\partial)$ which depends only on the
affine curve $X$ and not on the choice of its compactification,  
the projective curve $\lbar{X}$.

\end{prop}

The resulting homology groups of the chain complex 
$({\cal C}_\bullet(X),\partial)$ are denoted as before by $PH_k(X)$ and
are called polar homology groups of $X$ also in this case of an affine
$X$.

\smallskip


\begin{ex} \label{Z-P}
Let us consider a smooth projective curve of genus $g$ without a point,
$Z\smallsetminus\{P\}$. Then for the dimensions of polar homology
groups, $hp_k(X)=\dim HP_k(X)$, we get
$$
\begin{array}{lll}
hp_0(Z\smallsetminus\{P\}) = 2\,,   &
hp_1(Z\smallsetminus\{P\}) = 0\,,   & g=0\,,\\
hp_0(Z\smallsetminus\{P\}) = 1\,,   & 
hp_1(Z\smallsetminus\{P\}) = g-1\,, & g\geqslant 1\,.
\end{array}
$$
Indeed, the space $HP_1(Z\smallsetminus\{P\})$ is the space of
holomorphic 1-differentials on $Z$ which vanish at $P$. To calculate 
$HP_0(Z\smallsetminus\{P\})$  in the case $g\geqslant 1$ it is 
sufficient to notice that for any two points $Q_1,Q_2\in
Z\smallsetminus\{P\}$, the 0-cycle $(Q_1,q_1)+(Q_2,q_2)$ is
homologically equivalent to zero if and only if $q_1+q_2=0\,$
(the same condition as in the case of a non-punctured curve, 
cf. Example \ref{InCurv}). In the case
of $g=0$ an analogous statement requires three points to be involved
(unlike the case of a non-punctured projective line):
the corresponding 1-form on ${\C\P}^1$ has to have at least one 
zero, and hence at least three poles.
We collect the results about the curves in the following figures
(where we depict the complex curves by graphs, such that polar 
homology groups of the curves coincide with 
singular homology groups of the corresponding graphs).

\figone
\figtwo

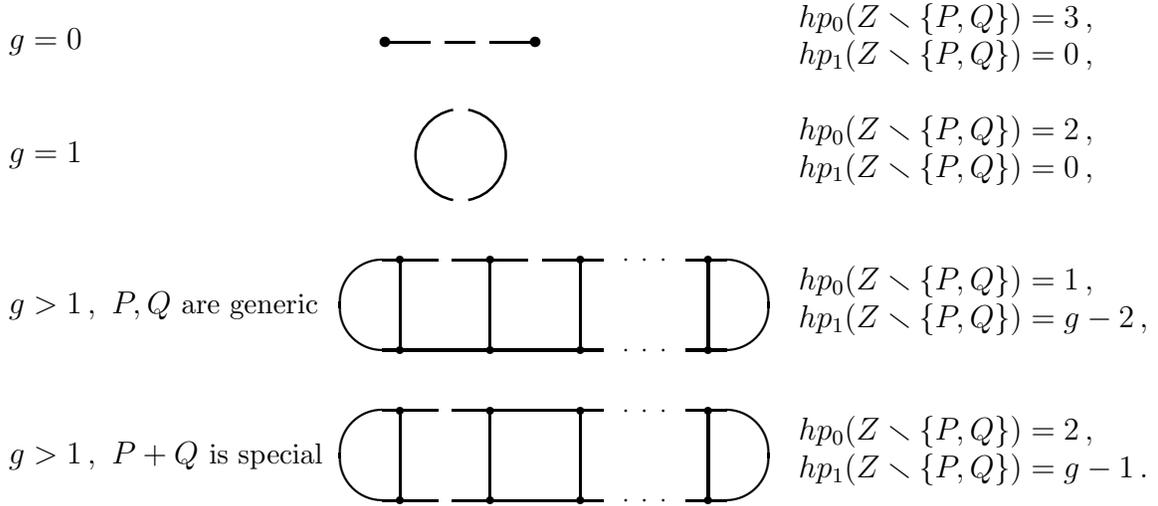
In a similar way, for a smooth projective curve without two points, 
$Z\smallsetminus\{P,Q\}$,
we get the results summarized in fig.~3. Here, one has
to distinguish the case of generic points $P$ and $Q$, and the case
when $P+Q$ is a special divisor and there are more 1-differentials with
zeros at $P,Q$ than generically. 

\figthree

\end{ex}


\begin{theo}     \label{MV}
{\bf (The Mayer--Vietoris sequence.)}
Let  a complex curve $X$ (either affine or projective) 
be the union of two Zariski open subsets $U_1$ and
$U_2\,$,  $X=U_1\cup U_2$ .
Then the following Mayer--Vietoris sequence 
of chains is exact:
$$
0 \to {\cal C}_k(U_1\cap U_2) \stackrel{\!\! i}{\longrightarrow}
{\cal C}_k(U_1)\oplus
{\cal C}_k(U_2) \stackrel{\!\!\sigma}{\longrightarrow}  {\cal C}_k(X)
\to 0.
$$
Here the map $\sigma$ represents  the sum of chains,
$$
\sigma:a\oplus b\mapsto a+b,
$$
and the map $i$ is the embedding of the chain lying in the intersection
$U_1\cap U_2$  as a chain in each subset $U_1$ and $U_2\,$:
$$
i: c \mapsto (c)\oplus(-c).
$$
This implies the following exact Mayer--Vietoris sequence in polar homology:
$$
\dots \to HP_k(U_1\cap U_2) \stackrel{\!\! i}{\longrightarrow}  HP_k(U_1)\oplus
HP_k(U_2) \stackrel{\!\!\sigma}{\longrightarrow}  HP_k(X)
\to HP_{k-1}(U_1\cap U_2)\to\dots
$$
\end{theo}

The proof of this theorem in the case of curves readily follows
from definitions of polar homology groups (by using a resolution
if the curve is singular). One can see that such a proof
essentially repeats the considerations with topological homology
of appropriate one-dimensional cell
complexes (i.e., graphs) as it is illustrated in Examples
\ref{InCurv} and \ref{Z-P}, as well as in figures 1-3 above.


\bigskip
\mysection{Connections and Gauge Transformations 
on Complex Curves and Surfaces} \label{Gauge}

The same philosophy of holomorphic orientation can be applied to 
field-theoretic notions in the following way. Suppose we have a functional
${\cal S}(\varphi) =\int_M L(\varphi, \partial_j\varphi)$ 
on {\it smooth} fields $\varphi$ (e.g., functions, connections, etc.)
on  a real (oriented) manifold $M$, and this functional 
is defined by an $n$-form $L$,
which depends on the fields and their derivatives.

Then on a complex $n$-dimensional
manifold $X$ equipped with a ``polar orientation'', i.e., with a
holomorphic or meromorphic $n$-form $\mu$, 
a complex counterpart ${\cal S}_\C$ of the functional 
$\cal S$ can be defined as follows:
${\cal S}_\C(\varphi) =\int_X\mu\wedge L(\varphi, {\bar\partial}_j\varphi)$.
Here $\varphi$ stands for {\it smooth } fields on a complex manifold $X$.
Now the $(0,n)$-form $L$ is integrated against the holomorphic orientation
$\mu$ over $X$.

Furthermore, the interrelation between the extremals of the real functional 
${\cal S}(\varphi)$ 
(on smooth  fields) on the real manifold $M$ and 
the boundary values of those fields
on  $\partial M$ (cf., e.g. \cite{S})
is replaced by the analogous interrelation
for the complex functional  ${\cal S}_\C(\varphi)$ (still on smooth fields)
on  $(X, \mu)$, i.e., a complex manifold $X$ equipped with 
polar orientation $\mu$,
and on its polar boundary, $\partial (X, \mu)
=( \text{div}_\infty\mu, \res~\mu)$.

Below we demonstrate some features of the above-mentioned 
parallelism for gauge transformations and   connections on curves 
and surfaces (cf. \cite{KR2, DT} for other examples).

\medskip

\begin{ssect}{Affine and double loop Lie algebras.} ~Our first example
is the correspondence between the affine Kac--Moody algebras on a circle 
($\R$-case) and the Etingof--Frenkel Lie algebras of currents over 
an elliptic  curve ($\C$-case) \cite{EF}.
\end{ssect}

We will use the following notations throughout this section.
Let $G$ be a simple simply connected 
Lie group that is supposed to be {\it compact} in the $\R$-case
and {\it complex} in the $\C$-case; 
${\mathfrak g}=\text{Lie}(G)$ its Lie algebra. 
Fix some smooth vector $G$-bundle $\cal E$ over an (either real or 
complex) manifold $M$. The notation $G^M$ (respectively,
${\mathfrak g}^M$) stands for the Lie group (respectively, Lie algebra) 
of $C^\infty$-smooth gauge transformations of $\cal E$.

\begin{defn} ~

$\R${\rm)} An affine Lie algebra ${\hat{\mathfrak g}^S}$ is the one-dimensional
central extension of the loop algebra 
${\mathfrak g}^S=C^\infty(S^1, \mathfrak g)$ (i.e., the gauge algebra 
over a circle) defined by the following 2-cocycle:
$$
c(U,V)=\int_{ S^1}tr(UdV) 
\qquad\text{for }~~ U,V\in {\mathfrak g}^S. 
$$ 

$\C${\rm)} {\rm\cite{EF}} An elliptic (or double loop) Lie algebra 
${\hat{\mathfrak g}}^{E}$ is a 
one-dimensional (complex) central extension of the gauge  algebra 
${\mathfrak g}^{E}$ over 
an elliptic curve $E$ by means of the following 2-cocycle:
$$
c(U,V)=\int_{E} \alpha\wedge tr(U{\bar\partial}V)
$$
where $\alpha$ is a holomorphic 1-form on $ E$ (its ``holomorphic 
orientation"), and  $U,V\in {\mathfrak g}^{ E}$.
\end{defn}

The original definition in [EF] was for the case of the current algebra
${\mathfrak g}^{ E}=C^\infty({ E}, \mathfrak g)$. However, it is  
valid in a more general case, which we need, for the group of gauge 
transformations of a bundle $\cal E$ not necessarily of degree zero.

The dual spaces to both affine and elliptic Lie algebras have a 
very natural geometric interpretation. Denote by ${\cal A}^M$
the infinite-dimensional affine space of all smooth connections 
(respectively, of all  (0,1)-connections) in the $G$-bundle $\cal E$ 
over real (respectively, complex) manifold $M$. 

Note that over a real curve   all connections are necessarily flat.
Analogously, over a complex curve every (0,1)-connection 
defines a structure of holomorphic bundle in $\cal E$
(since for such connections the curvature component 
$F^{0,2}$ is identically zero).
 
\begin{prop} ~

$\R${\rm)} (see,  e.g., {\rm\cite{PS}}) 
The space  ${\cal A}^S:=\{d+A~|~A\in\Omega^1(S^1,\mathfrak g)\}$ of 
 smooth $G$-connections over the circle $S^1$ 
can be regarded as (a hyperplane in) the dual space to the  
affine Lie algebra ${\hat{\mathfrak g}^S}$: the gauge transformations coincide
with the coadjoint action. 
Coadjoint orbits of the affine group, or the symplectic leaves of 
the linear Lie--Poisson
structure on the dual space $\left(\hat{\mathfrak g}^S\right)^*$,  
consist of gauge-equivalent 
connections and differ by (the conjugacy class of) the holonomy around $S^1$.

$\C${\rm)} {\rm\cite{EF}} The space of (0,1)-connections 
 $\{\bar\partial+A(z,\bar z)~|~A\in \Omega^{0,1}({E},\mathfrak g)\}$ 
in the bundle $\cal E$
over the elliptic curve $E$ can be regarded as  (a hyperplane in) 
the dual space $\left(\hat{\mathfrak g}^{E}\right)^*$ 
of the elliptic Lie algebra. 
The symplectic leaves of the Lie--Poisson structure 
in the dual space  $\left(\hat{\mathfrak g}^{E}\right)^*$
are enumerated by the equivalence classes of  
holomorphic $G$-bundles (or, different holomorphic structures in the 
smooth bundle $\cal E$) over the curve $E$.
\end{prop}

\begin{rem} Feigin and Odesski found a very interesting class 
of Poisson algebras (as well as their deformations, associative 
algebras) given by certain quadratic relations, and associated 
to a given complex $G$-bundle $\cal E$ over an elliptic curve $E$. It turned 
out that the symplectic leaves of those Poisson brackets are 
enumerated by  the isomorphism classes of holomorphic structures in $\cal E$ 
(see \cite{FO}), i.e., by the very same objects as the orbits of elliptic 
Lie algebras.  Therefore it would be interesting
to compare the transverse Poisson structures to the orbits 
of double loop Lie algebras with the transverse  structures
to the symplectic leaves of 
the Feigin--Odesski quadratic Poisson brackets. 
\end{rem}
 

\begin{ssect}{Gauge transformations over real surfaces.}
~Let $P$ be a real two-dimensional oriented manifold,
possibly with boundary $\partial P=\cup_j \Gamma_j$.
Let ${\cal A}^P$ be  the affine space of all smooth
connections in a trivial  $G$-bundle $\cal E$ over $P$.
It is convenient to fix any trivialization of $\cal E$ and identify
${\cal A}^P$ with the vector space $\Omega^1(P,\mathfrak g)$ of
smooth $\mathfrak g$-valued 1-forms on the surface:
$$
{\cal A}^P=\{ d+A~|~A\in\Omega^1(P,\mathfrak g)\}~.
$$
The space ${\cal A}^P$ is in a natural way a symplectic manifold
with the symplectic structure
$$
W:=\int_P\operatorname{tr}(\delta A\wedge \delta A)~,
$$
where $\delta$ is the exterior differential on ${\cal A}^P$, and
$\wedge$ stands to denote the wedge product both on ${\cal A}^P$
and $P$. The symplectic structure $W$ 
is invariant with respect to the gauge transformations
$$
A\mapsto g^{-1}Ag+g^{-1}dg~,
$$
where $g$ is an element of the group of gauge transformations,
$G^P$, i.e., it is a smooth $G$-valued function on the surface
$P$. However, this action is not Hamiltonian, if the surface $P$ 
has a non-empty boundary. In the latter case, the centrally extended group 
${\hat G}^P$ of gauge transformations on the surface  acts on ${\cal A}^P$
in a Hamiltonian way. 

We are interested in 
the quotient  of the  subset of flat connections
${\cal A}^P_{\text{fl}}\subset{\cal A}^P$ over the gauge groups action
of ${\hat G}^P$:
$$
{\cal M}^P_{\text{fl}}={\cal A}^P_{\text{fl}}/\hat{G}^P
=\{d+A\in {\cal A}^P\,|\,dA+A\wedge A=0\}/\hat{G}^P\,.
$$
The moduli space ${\cal M}^P_{\text{fl}}$ is a finite-dimensional  
manifold (with orbifold singularities); it can be also described
as the space of representations of the fundamental group $\pi_1(P)$ in $G$ modulo
conjugation.

The manifold ${\cal M}^P_{\text{fl}}$ can be endowed with a Poisson structure.
Its definition and properties can be conveniently dealt with by means of the 
Hamiltonian reduction ${\cal A}^P/\!/{\hat G}^P$. 
\end{ssect}

\begin{theo} ~

1) {\rm\cite{AB}} If the surface $P$ has no boundary, then the  space
${\cal M}^P_{\rm{fl}}$ of flat $G$-connections modulo  gauge transformations
on a surface $P$ is symplectic.

2) {\rm\cite{FR}} If $\partial P=\cup_j \Gamma_j$, then the moduli space 
${\cal M}^P_{\rm{fl}}$ 
on a surface $P$ with holes inherits a Poisson structure from
the space of all (smooth) $G$-connections.
The symplectic leaves of this structure are parameterised by the
conjugacy classes of holonomies around the holes (that is, a symplectic
leaf is singled out by fixing the conjugacy class of the
holonomy around each hole).
\end{theo}

We note that the second part of the theorem claims that the symplectic 
leaves of ${\cal M}^P_{\text{fl}}$ are labeled by 
the coadjoint orbits of the  affine Lie algebra 
on a circle (or of several copies of the
affine algebra, with each copy situated at a different boundary component 
of the surface $P$), 
since those orbits are parameterised by the conjugacy classes of holonomies
around the circle.


\bigskip
\begin{ssect}{Gauge transformations over complex surfaces.}
In this section we present a complex counterpart of the description 
of the Poisson structures on moduli spaces. Let $Y$ be 
a compact {\it complex} surface ($\dim_{\C}\,Y=2$).
Choose a polar analogue of orientation, i.e., a holomorphic or meromorphic
2-form $\beta$ on $Y$.  Let $\beta$ be a meromorphic 2-form on $Y$, which has
only first order poles on a smooth curve $X$.
The curve $X\subset Y$ will play the
role of the boundary of the surface $Y$ in our considerations.  Moreover, 
assume  that $\beta$ has no zeros (the situation analogous to a
smooth oriented compact real surface).  
Then $X$ is an anticanonical divisor in $Y$ and it
has to be an elliptic curve $E$, or, may be, a number of non-intersecting
elliptic
curves. (Example: $Y={\C\P}^2$ with a smooth cubic 
as an anticanonical divisor. As a matter of  fact, many Fano surfaces fall 
into this class. )
If it happens that $\beta$ has no zeros and no poles (i.e., $Y$
is ``oriented, without boundary") it means that we deal with either a K3
or an abelian surface. 
\end{ssect}

Let $\cal E$ be a smooth vector $G$-bundle over $Y$ which can be endowed with
a holomorphic structure and $\text{End}\,\cal E$ be the corresponding bundle of
endomorphisms with the fiber ${\mathfrak g}=\text{Lie}(G)$. Let ${\cal A}^Y$
denote the infinite-dimensional affine
space of smooth $(0,1)$-connections in $\cal E$. By choosing a
reference holomorphic structure $\bar\partial_{\,0}$,
$\bar\partial_{\,0}^{\;2}=0$, in
$\cal E$, the space ${\cal A}^Y$ can be identified with the vector space
$\Omega^{(0,1)}(Y,\text{End}\,\cal E)$ of $(\text{End}\,\cal E)$-valued
$(0,1)$-forms on
$Y$, i.e.,
$$
{\cal A}^Y = \{\bar\partial_{\,0} +A\,|\, A\in
\Omega^{(0,1)}(Y,\text{End}\,{\cal E})\}~.
$$
We often shall write $\bar\partial$ instead of $\bar\partial_{\,0}$,
keeping in mind that this corresponds to a reference holomorphic
structure in $\cal E$ when it applies to sections of $\cal E$ or associated
bundles.

The space ${\cal A}^Y$ possesses a natural holomorphic symplectic
structure
$$
 W_{\C} := \int_Y \beta \wedge\operatorname{tr}(\delta A_1\wedge
\delta
A_2)~,
$$
where $\beta$ is the ``polar orientation" of $Y$, while the other
notations are the same as above.
The symplectic structure $W_{\C}$ is invariant with
respect to the gauge transformations
$$
A\mapsto g^{-1}Ag+g^{-1}\bar\partial g~,
$$
where $g$ is an element of the group of gauge transformations,
i.e., the group of automorphisms of the smooth bundle $\cal E$. Abusing
notation we denote this group by $G^Y$.

Again, we will need to centrally extend the group $G^Y$ of gauge 
transformations to make the action Hamiltonian.
In the momentum map, taking the curvature is replaced by
the mapping: $A\mapsto \beta\wedge F^{0,2}(A)=\beta\wedge 
(\bar\partial A+A\wedge A)$.
When equating the result to zero, instead of the flatness condition
$F(A)=0$,  we come to the relation $F^{0,2}(A)=0$, which 
singles out (0,1)-connections 
defining {\it holomorphic} structures in $\cal E$. Denote 
the space of such $\bar\partial$-connections by ${\cal A}^Y_{\text{hol}}$. 
The set of isomorphism classes of holomorphic structures in $\cal E$
is represented by the quotient 
$$
{\cal A}^Y_{\text{hol}}/{\hat G}^{Y}=
\{\bar\partial +A\in {\cal A}^Y\,|\,\bar\partial A+A\wedge A=0\}/
{\hat G}^{Y}\,.
$$

Analogously to the moduli space of flat connections on a real surface,
we would like to study the Poisson geometry of the moduli space of  
holomorphic bundles over a complex surface. However, the question 
of existence and 
singularities of such a moduli space is much more subtle.
Suppose the bundle $\cal E$ was chosen in such a way that there exists
some version of the moduli space of holomorphic structures in $\cal E$
(e.g., (semi-)stable bundles). Denote by ${\cal M}^Y_{\text{hol}}$
the non-singular part of that moduli space. This finite-dimensional manifold
can be equipped with a holomorphic Poisson structure.

Since ${\cal M}^Y_{\text{hol}}$ is an open dense subset in the space of isomorphism classes of 
holomorphic bundles,
$$
{\cal M}^Y_{\text{hol}}\subset {\cal A}^Y_{\text{hol}}/{\hat G}^{Y}\,,
$$
the  Poisson structure on ${\cal M}^Y_{\text{hol}}$ 
can be studied by means of the Hamiltonian reduction.

\begin{theo}{\label{holPS}} ~

1) {\rm\cite{Mu}} If $Y$ is a K3 surface or a complex torus of dimension 2,
i.e., if the 2-form $\beta$ is holomorphic on $Y$, then the moduli space
${\cal M}^Y_{\rm{hol}}$ admits a holomorphic symplectic structure.

2) If $\beta$ is meromorphic, the  moduli 
space ${\cal M}^Y_{\rm{hol}}$ of
holomorphic bundles possesses a (holomorphic) Poisson structure (see
{\rm\cite{Bon, Bot, Tu}}, where
the Poisson structure is given in intrinsic terms).  The
symplectic leaves of this structure are parameterised by the 
isomorphism classes of the restrictions of  bundles to 
the anticanonical divisor $X\subset Y$, {\rm\cite{KR2}}.
\end{theo}

 Thus, the symplectic leaves of the Poisson structure on 
${\cal M}^Y_{\text{hol}}$
are distinguished by the moduli of holomorphic bundles on  elliptic curve(s)
$X$, or, which is the same, by 
coadjoint orbits of the corresponding elliptic algebras
${\hat{\mathfrak g}}^{X}$ on (the connected components of) the smooth
divisor $X\subset Y$.
\footnote{Note that the choice of isomorphism classes of bundles on $X$  
must be subject to the condition that they
arise as restrictions of bundles defined over $Y$.}

The above consideration can be
extended with minimal changes to the case of a non-smooth divisor $X$,
in particular, to $X$ consisting of several
components intersecting transversally. 
(Example: $Y={\C\P}^2$ with $\beta=dxdy/xy$.) 
In the latter case the corresponding degeneration of the
elliptic algebra ${\hat{\mathfrak g}}^{X}$ can be described in terms of
(several copies of) the  current algebra on a punctured 
${\C\P}^1$.
\medskip


\begin{ssect}{Chern--Simons functionals.}
First, let $M$ be a real compact three-dimensional manifold
with boundary $P=\partial M$, and $\cal E$ a trivial $G$-bundle over $M$. 
The {\it Chern--Simons functional} CS on 
the space of $G$-connections ${\cal A}^M$ is given by the formula
$$
\text{CS}(A)=\int\limits_M {tr}(A\wedge dA+\frac{2}{3}A\wedge A\wedge A).
$$
Extremals of this functional are flat connections on $M$. 
An action functional on the  fields in three dimensions 
defines a symplectic structure on the space of fields in two
dimensions. 
(It arises due to the relation between the boundary values of the
fields and the solutions to the Euler-Lagrange equations;
this is essentially the Hamiltonian approach to the
corresponding variational problem.) 
In the present case, as it is well known \cite{W2}, the corresponding
symplectic manifold is the moduli space of flat connections ${\cal
M}^P_{\rm fl}$ on $P=\partial M$.

The path integral corresponding  to the Chern-Simons functional
can be related with invariants of links in a three-dimensional manifold
\cite{W2}
and, in the simplest case of an abelian gauge group  
($G=U(1)$), 
reproduces the definition 
of the Gauss linking number.
\end{ssect}

The ``holomorphic'' counterpart of the Chern--Simons
functional,
$$
\text{CS}_{\C}(A)=\int\limits_Z\gamma\wedge {tr}(A\wedge\bar\partial A
+\frac{2}{3}A\wedge A\wedge A)
$$
suggested in ref.~\cite{W}, can be treated to some extent similarly. 
Here $\text{CS}_{\C}(A)$ is considered as a functional on the space of 
$\bar\partial$-connections $A\in {\cal A}^Z$ in a trivial $G$-bundle 
$\cal E$ over 
a complex 3-fold $Z$, where $Z$ is  equipped with 
a meromorphic (``polar orientation'') 3-form $\gamma$ without zeros,
but, may be, with poles of the first order. 
In such a situation, one can apply the arguments similar 
to the case of the ordinary Chern--Simons theory, provided that one
replaces everywhere $d$ by $\bar\partial$ and, instead of real boundary, 
one deals with the polar boundary $Y:= \dvsr_\infty\gamma\subset Z$.
The extrema of $\text{CS}_{\C}(A)$ are given now by integrable 
$\bar\partial$-connections ($\bar\partial_A^2=0$), that is by
holomorphic bundles over $Z$ (which are counterparts 
of flat connections in three dimensions).
Then, at the complex two dimensional ``boundary'', one gets 
the symplectic manifold 
${\cal M}^Y_{\text{hol}}$ of moduli of holomorphic bundles over a
complex surface $Y$ (as a counterpart of the moduli space of flat
connections in two real dimensions).

The holomorphic Chern--Simons theory in the case of an abelian gauge
group $G$ on a complex simply connected 3-fold $Z$ 
can be discussed even further, at the level of path integrals, 
without much difference with its ``real''
prototype (unlike the case of an arbitrary, non-abelian, gauge group
$G$, which is much more complicated and still lacks a rigorous treatment),
cf. \cite{FT, T}.

\medskip

\begin{ssect}{Polar links.} 
In the abelian case, the quantum holomorphic Chern--Simons theory 
reproduces a holomorphic analogue of the linking number. Its
definition can be immediately found, again, by analogy with the ordinary one.

Let $Z$ be a complex projective three-dimensional manifold,
equipped, as above, with a meromorphic 3-form $\gamma$ without zeros.
Consider two smooth polar 1-cycles $(C_1,\alpha_1)$ and $(C_2,\alpha_2)$
in $Z$, i.e., $C_1$ and $C_2$ are smooth complex curves equipped with
holomorphic 1-forms.
Let us take the 1-cycles which are polar boundaries. This means, in
particular, that there exists such a 2-chain $(S_2,\beta_2)$  that
$(C_2,\alpha_2)=\partial\,(S_2,\beta_2)$.
Suppose, the curves $C_1$ and $C_2$ have no common points and
$S_2$ is a smooth surface which
intersects transversely with the curve $C_1$. Then, we define
the {\it polar linking} number of the 1-cycles above
as the polar intersection number (cf.\ (\ref{X})) 
of the 2-chain  $(S_2,\beta_2)$
with the 1-cycle $(C_1,\alpha_1)$:
$$
\lk\left( (C_1,\alpha_1), (C_2,\alpha_2)\right):=
\sum\limits_{P\in C_1\cap \,S_2} 
\frac{\alpha_1(P)\wedge\beta_2(P)}{\gamma(P)}\,.
$$
One can show that the expression above does not depend on the choice of
$(S_2,\beta_2)$,  and has certain invariance
properties mimicking those of the topological linking number
within  the framework of the ``polar'' approach.
\end{ssect}


\bigskip

\begin{ackn}
B.K. and A.R. are grateful for hospitality to the Max-Planck-Institut f\"ur
Mathematik in Bonn, where a large part of this work was completed.
The present work was partially sponsored by PREA and McLean awards.

We are grateful to B.~Feigin, R.~Thomas, and A.~Odesski
for interesting discussions. 

The work of B.K. was partially supported by an Alfred P. Sloan
Research Fellowship,   by the NSF and NSERC research grants.
The work of A.R. was supported in part by the Grants RFBR-00-02-16530,
INTAS-99-1782 and the Grant 00-15-96557 for the support of scientific
schools.

\end{ackn}

\begin{smallbibl}{MMMM}

\bibitem[A]{Arn}\ V.I.~Arnold,
{\em  Arrangement of ovals of real plane algebraic curves,
 involutions
of smooth four-dimensional  manifolds, and on arithmetic
of integral-valued quadratic forms}\/,
Func. Anal. and Appl.  {\bf 5:3}  (1971),  169--176

\bibitem[AB]{AB}\ M.\,Atiyah and R.\,Bott
{\em The Yang-Mills equations over a Riemann surface}\/, 
Philos. Trans. Roy. Soc. London {\bf 308} (1982) 523--615

\bibitem[Bon]{Bon}\ A.\,Bondal, {\em unpublished}

\bibitem[Bot]{Bot}\ F.\,Bottacin,
{\em Poisson structures on moduli spaces of sheaves over Poisson surfaces}\/, 
Invent. Math. {\bf 121:2} (1995) 421--436

\bibitem[DT]{DT}\ S.K.~Donaldson and R.P.~Thomas,
{\em Gauge theory in higher
dimensions}\/,  The geometric universe (Oxford, 1996),
Oxford Univ. Press, Oxford (1998) 31--47.

\bibitem[EF]{EF}\ P.I.\,Etingof and I.B.\,Frenkel,
{\em Central extensions of current groups in two dimensions}\/, 
Commun. Math. Phys. {\bf 165} (1994) 429-444

\bibitem[FO]{FO}\ B.L.~Feigin and A.V.~Odesski
{\em Vector bundles on an elliptic curve and Sklyanin algebras}\/, 
Amer. Math. Soc. Transl. Ser. 2 {\bf 185} (1998) 65--84 ~(q-alg/9509021)

\bibitem[FR]{FR}\ V.V.\,Fock and A.A.\,Rosly,
{\em Poisson structures on moduli of flat connections on
Riemann surfaces and $r$-matrices}\/, 
preprint (1992); {\em Flat connections and polyubles}
Theor. Math. Phys. {\bf 95:2} (1993) 526--534

\bibitem[FK]{FK}\ I.B.~Frenkel and B.A.~Khesin,
{\em Four-dimensional realization of two-dimensional current
groups}\/, Commun. Math. Phys. {\bf 178} (1996) 541--561

\bibitem[FT]{FT}\ I.B.~Frenkel and A.N.~Todorov,~
{\em paper in preparation}

\bibitem[Gr]{Gr}\ P.A.~Griffiths,~ {\em Variations on a Theorem of
Abel}\/, Invent. Math., {\bf 35} (1976) 321--390

\bibitem[K]{Kh}\ B.A.\,Khesin,~ 
{\em Informal complexification and Poisson structures on moduli spaces}\/, 
AMS Transl., Ser. 2, {\bf 180} (1997) 147--155

\bibitem[KR]{KR}\ B.~Khesin and A.~Rosly, {\em Symplectic geometry on
moduli spaces of holomorphic bundles over complex surfaces}\/,
The Arnoldfest,  Toronto 1997, Editors: E.~Bierstone et al.,
Fields Institute/AMS Communications (1999) {\bf 24}, 311--323.

\bibitem[KR2]{KR2}\ B.~Khesin and A.~Rosly,~ 
{\em Polar homology}\/,
Preprint, math.AG/0009015 (2000)

\bibitem[LMNS]{arhar}\ A.~Losev, G.~Moore, N.~Nekrasov, and S.~Shatashvili,
{\em Four--Dimensional Avatars of Two-Dimensional RCFT}\/,
Nucl. Phys. Proc. Suppl. {\bf 46} (1996) 130-145
~(hep-th/9509151)

\bibitem[Mu]{Mu}\ S.\,Mukai,
{\em Symplectic structure of the moduli space of stable sheaves on
    an abelian or $K3$ surface}\/, 
Inven. Math. {\bf 77} (1984) 101--116

\bibitem[PS]{PS} A.~Pressley and G.~Segal,
{\em Loop groups}\/,
 Clarendon Press, Oxford (1986)

\bibitem[S]{S}\ A.S.~Schwarz, 
{\em Symplectic formalism in conformal field theory}\/, Sym\'etries
quantiques (Les Houches, 1995), North-Holland, Amsterdam (1998) 957--977

\bibitem[T]{T}\ R.P.~Thomas,
{\em Gauge theory on Calabi--Yau manifolds}\/, Ph.D. thesis, Oxford
(1997), 1--104

\bibitem[Tyu]{Tu} A.\,Tyurin,
{\em Symplectic structures on the moduli spaces of vector bundles
on algebraic surfaces with $p_g>0$}\/,
Math. Izvestia {\bf 33} (1987) 139-177

\bibitem[W1]{W2}\ E.~Witten,
{\em Quantum field theory and the Jones polynomial}\/,
Commun. Math. Phys., {\bf 121:3} (1989) 351--399

\bibitem[W2]{W}\ E.~Witten,
{\em Chern--Simons gauge theory as a string theory}\/,
The Floer memorial volume, Progr. Math., {\bf 133}, Birkh\"auser,
Basel (1995) 637--678 ~(hep-th/9207094)

\end{smallbibl}

\end{document}